\documentclass{amsart}
\usepackage{psfig}

\newfam\Bbbfam
\def\R{{\mathbb R}}
\def\grad{\nabla}		  

\makeatletter
\let\ii\item
\newtheorem{crlry}{Corollary}[section]
\def\statecorollary{\@ifnextchar[{\@statecorollary}{\nr@statecorollary}}
\long\def\@statecorollary[#1]#2{\begin{crlry}\label{#1}#2\end{crlry}}
\long\def\nr@statecorollary#1{\begin{crlry}#1\end{crlry}}
\def\statecorollarypf{\@ifnextchar[{\@statecorollarypf}{\nr@statecorollarypf}}
\long\def\@statecorollarypf[#1]#2{\begin{crlry}\label{#1}#2\end{crlry}\proof}
\long\def\nr@statecorollarypf#1{\begin{crlry}#1\end{crlry}\proof}

\newtheorem{thm}{Theorem}[section]
\def\statetheorem{\@ifnextchar[{\@statetheorem}{\nr@statetheorem}}
\long\def\@statetheorem[#1]#2{\begin{thm}\label{#1}#2\end{thm}}
\long\def\nr@statetheorem#1{\begin{thm}#1\end{thm}}
\def\statetheorempf{\@ifnextchar[{\@statetheorempf}{\nr@statetheorempf}}
\long\def\@statetheorempf[#1]#2{\begin{thm}\label{#1}#2\end{thm}\proof}
\long\def\nr@statetheorempf#1{\begin{thm}#1\end{thm}\proof}

\newtheorem{lmma}{Lemma}[section]
\def\statelemma{\@ifnextchar[{\@statelemma}{\nr@statelemma}}
\long\def\@statelemma[#1]#2{\begin{lmma}\label{#1}#2\end{lmma}}
\long\def\nr@statelemma#1{\begin{lmma}#1\end{lmma}}
\def\statelemmapf{\@ifnextchar[{\@statelemmapf}{\nr@statelemmapf}}
\long\def\@statelemmapf[#1]#2{\begin{lmma}\label{#1}#2\end{lmma}\proof}
\long\def\nr@statelemmapf#1{\begin{lmma}#1\end{lmma}\proof}

\def\putlabel(#1,#2)[#3]#4{\put(#1,#2){\makebox(0,0)[#3]{#4}}}

\def\@loadps#1{\@ifnextchar[{\arg@loadps{#1}}{\do@loadps{#1}{\@aspicwd}{\@aspicht}}}
\def\arg@loadps#1[#2,#3]{\do@loadps{#1}{#2\unitlength}{#3\unitlength}}
\def\do@loadps#1#2#3{\makebox(0,0)[bl]{\psfig{figure=#1,width=#2,height=#3}}}
\def\loadps#1{\typeout{LOADPS: loading #1}\@loadps{#1}}

\makeatother

\newcommand{\q}{{\it\bf q}}
\newcommand{\p}{{\it\bf p}}

\newcommand{\s}{{\it\bf s}}
\newcommand{\x}{{\it\bf x}}

\newcommand{\bu}{{\it\bf u}}
\newcommand{\ba}{{\it\bf a}}
\newcommand{\bbb}{{\it\bf b}}
\newcommand{\bC}{{\it\bf M}}
\newcommand{\lc}{\ell_{ \bC}}
\newcommand{\rc}{r_{ \bC}}
\newcommand{\li}{\ell_{\cap}}
\newcommand{\ri}{r_{\cap}}
\newcommand{\cd}{{{}_\star}}

\title{Symbolic Dynamics of the Collinear Three-Body
Problem}
\author{Samuel R. Kaplan}

\address{Department of Mathematics, Bowdoin College\footnote{The
author is presently at the University of North Carolina at Asheville,
Asheville, NC 28801} Brunswick, ME 04011}

\begin{document}

\begin{abstract}
Solutions to the collinear three-body problem which do not end in
triple collision pass through an infinite number of binary
collisions. Given three masses, we show that four geometric quantities
generate a finite description of itineraries of binary collisions. In
the best circumstances, this description is semi-conjugate to a
Poincar\'{e} map of the flow. For other cases these quantities give
upper and lower bounds on the itineraries which can occur. In addition
to describing the dynamics of the collinear three-body problem, the
results of this paper rederives the existence of oscillatory motion in
the $N$-body problem for $N\geq 3$.
\end{abstract}

\maketitle

\section{Introduction}In the seventeenth century, Newton formulated the universal law of
gravitation and completely solved the two-body problem. Moreover, his
methods confirmed Kepler's laws of planetary motion. After this grand
success Newton turned his fluxions to the motions of the Earth, moon
and sun. He eventually gave up working on this three-body problem,
saying it gave him headaches \protect\cite{MH}.

In the intervening years much work was done on the three-body problem
generating new techniques and new questions.  After
Poincar\'{e} showed that chaotic behavior can occur in the restricted
three-body problem \protect\cite{Po}, the search for a complete solution halted
abruptly. Poincar\'{e}'s result changed the issue from trying to solve
an initial value problem to asking what behaviors can occur in the
three-body problem.

A critical step in understanding what dynamics can occur in the
three-body problem was identifying behavior near
triple-collision. McGehee made a great leap forward in this field in
1974 with his analysis of triple-collision behavior in the collinear
three-body problem \protect\cite{Mc}. He introduced a change of
variable to what is now known as McGehee coordinates, in which the
differential equations for Newton's universal law of gravitation are
extended to triple collision. By understanding the dynamics at triple
collision, one understands the dynamics near triple collision via
continuity of the flow.

Generating a global analysis of the collinear three-body problem
requires an understanding of how near-triple collision behavior
effects the entire flow. Using ideas of Saari and Xia
\protect\cite{SX}, Meyer and Wang undertook this endeavor
\protect\cite{MW}, giving a nice picture of the geometry of the phase
space for the collinear three-body problem.

Meyer and Wang's method of analysis involves generating a Poincar\'{e}
slice to the flow and partitioning that slice with pieces of stable
manifold for triple collision. Then in at least some part of the flow
the regions generate a subshift of finite type on an infinite set of
symbols. This is enough to show that chaotic behavior occurs in the
collinear three-body problem. However, this result raises the question
of determining exactly what itineraries of binary collisions are
allowed and how the set of allowed itineraries changes as the masses
vary.

In order to gain insight on the role of geometry in the collinear
three-body problem, the author chose to use McGehee coordinates. In
contrast, Meyer and Wang choose coordinates so that binary collision
is represented by a sink. Their approach has some advantages; however,
in McGehee coordinates binary collisions are represented by
half-planes. Thus, McGehee coordinates reveal more detail on how the
stable manifold for triple collision intersects the Poincar\'{e}
slice.

The crucial construction in this paper is the partition of a
Poincar\'{e} slice into a finite number of regions bounded by pieces
of stable manifold for triple collision. The Theorems in Section~5
give conditions on when this partition is Markov. When the partition
is Markov, we can construct a semi-conjugacy from the Poincar\'e map
to a directed graph and exactly describe the set of allowed
itineraries. When the partition is not Markov, we can construct two
sub-shifts of finite type, one containing all allowed sequences and
one containing guaranteed sequences. These two sub-shifts serve as
upper and lower bounds on the set of allowed itineraries.

For details of this approach applied to a model problem in which more
detail can be explicitly computed, see \protect\cite{K1, K2, K3, K4} which
consider the dynamics of the collinear one-bumper two-body problem. 

Except for very special cases, the conditions required by the theorems
below can only be computed numerically. Even so, the theorems in this
paper give a picture of what symbolic dynamics occur in the collinear
three-body problem and what geometric quantities one needs to compute
in order to describe those dynamics.

The symbolic dynamics of the collinear three-body problem are enough
to guarantee the existence of special solutions (oscillatory motion)
in the $N$-body problem (see the Corollary in Section~6).

\section{Hamiltonian Coordinates}We begin by stating the collinear three-body problem in Hamiltonian
coordinates and then look at the regularized system in McGehee coordinates. 

Three points masses have masses $m_i>0$ and positions
$q_1\leq q_2\leq q_3\in\R$. The potential energy is given by  
\begin{eqnarray}
U&=&\sum_{i>j} \frac{m_i m_j}{q_i-q_j} \protect\label{eq:pot}.
\end{eqnarray}

The motion of the particles under gravitational force is described by
the systems of differential equations
\begin{eqnarray}
m_i \ddot{q_i}=\nabla_{q_i} U, \qquad\qquad i=1, 2, 3 \protect\label{eq:ham}
\end{eqnarray}
where $\grad_{q_i} U$ is the gradient of $U$ with respect to $q_i$. 

A position, $(q_1,q_2,q_3)$ is called a {\it binary collision} if
either $q_1=q_2$ or $q_2=q_3$. If $q_1=q_2=q_3$ the position is called
a {\it triple collision}. The above system is defined everywhere
except at binary and triple collisions. Given an initial positions
(not at collision) and momenta at time $t=0$, a unique solution exists
on a maximal interval $[0,t^*)$. If $t^*<\infty$ then the solution is
said to have a singularity at $t^*$. The only singularities which
occur in the collinear three-body problem are due to collision \protect\cite{Pa},
though non-collision singularities have been found for other $N$-body
problems \protect\cite{X}.

Double collisions can be regularized. That is, one can change
variables so that double collision transforms to a regular point of
the flow \protect\cite{E}. This extension corresponds to an elastic bounce. McGehee
further showed that triple collision is not regualrizable. However, in
McGehee coordinates, the flow is bounded by an invariant compact
manifold which correspond to triple collision. The flow on this
invariant manifold, called the collision manifold, guides solutions
which pass near triple collision. 

We change from Hamiltonian to McGehee coordinates via three
steps. First is a change of variables to polar coordinates. Second is
a time scaling change of variable so that solutions slow down as they
pass near triple collision. Finally, one regularizes binary
collisions.

The remainder of this section is a review of the change of variables
to McGehee coordinates. Readers familiar with this material should
feel free to pass on to the next section.

Let $\q=(q_1,q_2,q_2)\in\R^3$ be the vector of positions. Define
$p_i=m_i \dot{q_i}$ and let $\p=(p_1,p_2,p_3)\in\R^3$ be the vector of
momenta. Let 
$$
M=\left(\begin{array}{ccc} 
	m_1 & 0 & 0\\
	0 & m_2 & 0\\
	0 & 0 & m_3 
	\end{array}
   \right).
$$
Then Equations~\protect\ref{eq:pot} and~\protect\ref{eq:ham} can be rewritten as
\begin{eqnarray}
U(\q) &=& \sum_{i>j} \frac{m_i m_j}{q_i-q_j} \nonumber\\[.15in]
M \ddot{\q} &=& \grad U(\q). \protect\label{eq:grad}
\end{eqnarray}

We can also write the kinetic energy for the system as 
$$
T(\p)=\makebox{$\frac12$} \p^T M \p.
$$

The Hamiltonian for the system is 
\begin{eqnarray*}
H(\q,\p)=T(\p)-U(\q),
\end{eqnarray*}
and Equation~\protect\ref{eq:grad} can be written as the system
\begin{eqnarray}
\dot{\q}&=&H_{\p}(\q,\p) \:=\: M^{-1}\p \nonumber\\
\dot{\p}&=&-H_{\q}(\q,\p) \:=\: \grad U(\q). \protect\label{eq:linsys}
\end{eqnarray}

The function $T$ is defined everywhere in $\R^3$. The function $U$ is
defined everywhere except at collisions. 

Next we break up $\q$ into radial and angular components. Define
$$r=\left(\q^TM\q\right)^{1/2}.$$ Notice that a level set of $r$ is an
ellipsoid in $\R^3$. Let $S=\{\q\:|\:r=1\}$. Then a point in $S$ is called
a configuration for the system of particles.

We now define the variables:
\begin{eqnarray*}
r&=&\left(\q^T M \q\right)^{1/2}\\
\s&=& r^{-1} \q\\
y&=& \p^T s\\
\x&=& \p-yM\s.
\end{eqnarray*}
In these coordinates, $r$ is the size of the configuration in a
inertial norm, $\s\in S$ is the configuration and represents the
direction of $\q$, $y$ is the projection of $\p$ in the direction of
$\q$ and $\x$ represents the direction of change in the
configuration. Note that $\s$ and $\x$ are orthogonal. 

In the new polar coordinates, the Hamiltonian, $H(\q,\p)=h$ can be
written as
$$
\makebox{$\frac12$} (\x^T M^{-1} \x +y^2)-r^{-1} U(\s)=h
$$
and the equations of motion become
\begin{eqnarray*}
\dot{r}&=&y  \\
\dot{y}&=&r^{-1} \x^T M^{-1} \x-r^{-2} U(\s)\\
\dot{\s}&=&r^{-1} M^{-1} \x\\
\dot{\x}&=&-r^{-1}y\x-r^{-1}(\x^T M^{-1}\x)+r^{-2}U(\s)M\s+r^{-2}\grad
U(\s). 
\end{eqnarray*}

We next introduce two new coordinates, $\bu=r^{1/2} \x$ and $v=r^{1/2}
y$ and scale time via $dt = r^{3/2} dt'$. In these time-scaled
coordinates, the Hamiltonian, $H(\q,\p)=h$ can be written as
$$
\makebox{$\frac12$} (\bu^T M^{-1} \bu +v^2)- U(\s)=r h
$$
and the equations of motion become
\begin{eqnarray*}
\dot{r}&=&r v  \\
\dot{v}&=&\makebox{$\frac12$} v^2+\bu^T M^{-1} \bu- U(\s) y\\
\dot{\s}&=&M^{-1} \bu\\
\dot{\bu}&=&-\makebox{$\frac12$} v\bu-(\bu^T M^{-1}\bu)M\s+U(\s)M\s+\grad U(\s). 
\end{eqnarray*}

We reduce the dimension of System~\protect\ref{eq:linsys} by fixing the center
of mass at the origin and setting the total momentum to zero. This
reduces the problem to a four-dimensional phase space. 

To express this reduced system, note that there is a unique point on
$S$ so that $q_1=q_2<q_3$ and the center of mass, $M\q=0$. Call this
unique point $\ba=(a_2,a_2,a_3)$. Likewise, there is a unique point on
$S$ so that $q_1<q_2=q_3$ and the center of mass, $M\q=0$. Call this
unique point $\bbb=(b_1,b_2,b_2)$. One can compute that $0<\ba^T
M\bbb <1$. Let $\lambda$ be the least positive solution to 
$$
\cos(2\lambda)=\ba^T M \bbb .
$$

We now introduce an angular potential function, 
\begin{eqnarray*}
W(s)&=&\frac{2}{\lambda}(W_1(s)+W_2(s)+W_3(s))\sin(2 \lambda) 
\end{eqnarray*}
where
\begin{eqnarray*}
W_1(s)&=& \frac{m_1 m_2 (1-s)}{(b_2-b_1) Sn(\lambda(1+s))}\\
W_2(s)&=& \frac{m_2 m_3 (1+s)}{(a_3-a_2) Sn(\lambda(1-s))}\\
W_1(s)&=& \frac{\lambda m_1 m_3 (1-s^2)}{(b_2-b_1)
\sin(\lambda(1+s))+(a_3-a_2)\sin(\lambda (1-s))}
\end{eqnarray*}
and 
$$
Sn(x)=\frac{\sin(x)}{x}.
$$

To regularize double collisions, we introduce a new variable, 
$$
w=\frac{(1-s^2)u}{\sqrt{W(s)}}.
$$
and scale time again via, 
$$
dt'=\frac{\lambda (1-s^2)}{\sqrt{W(s)}} d\tau.
$$

Finally, we are ready to write System~\protect\ref{eq:linsys} in McGehee coordinates.
\begin{eqnarray}
\frac{dr}{d\tau} &=& \frac{\lambda (1-s^2)}{\sqrt{W(s)}} r v
\nonumber\\
\frac{dv}{d\tau} &=& \frac{\lambda}{2} \sqrt{W(s)}
\left(1-\frac{(1-s^2)}{W(s)}(v^2-4rh)\right) \nonumber\\
\frac{ds}{d\tau} &=& w \nonumber \\
\frac{dw}{d\tau} &=&
-s+\frac{2s(1-s^2)}{W(s)}(v^2-2rh)+\frac{W'(s)}{2W(s)}(1-s^2-w^2)-\frac{\lambda
(1-s^2)}{2\sqrt{W(s)}}vw. \protect\label{eq:mcg}
\end{eqnarray}
In McGehee coordinate, the Hamiltonian, $H(\q,\p)=h$ can be rewritten
as
\begin{eqnarray}
W(s)(w^2+s^2-1)+(1-s^2)^2(v^2-2rh)=0 \protect\label{eq:w}
\end{eqnarray}

System~\protect\ref{eq:mcg} is defined for all values of $r, v, s$ and $w$. Note
that $s$ is the configuration coordinate and varies from $-1$ to
$1$. The configuration $s=-1$ corresponds to $q_1=q_2<q_3$, a left
binary collision. The configuration $s=1$ corresponds to
$q_1<q_2=q_3$, a right binary collision. Moreover, the vector fields
is now defined at $r=0$, triple collision. 

When $r=0$ the Hamiltonian $H=h$ in Equation~\protect\ref{eq:w} yields the relation
\begin{eqnarray}
W(s)(w^2+s^2-1)+(1-s^2)^2(v^2)=0. \protect\label{eq:energy}
\end{eqnarray}
Equation~\protect\ref{eq:energy} defines a manifold, $\bC$, in $\R^3$, called
the {\it collision manifold}, which is independent of the total energy,
$h$. The phase space for System~\protect\ref{eq:mcg} is bounded by $\bC$.

\section{McGehee Coordinates}We briefly describe the flow in McGehee coordinates, summarizing
results from McGehee's research \protect\cite{Mc}. For the rest of the paper, we assume
that the total energy, $h$, is negative.

\begin{floatfig}
\begin{picture}(3.94,2.42)
\put(0,-2.42){\loadps{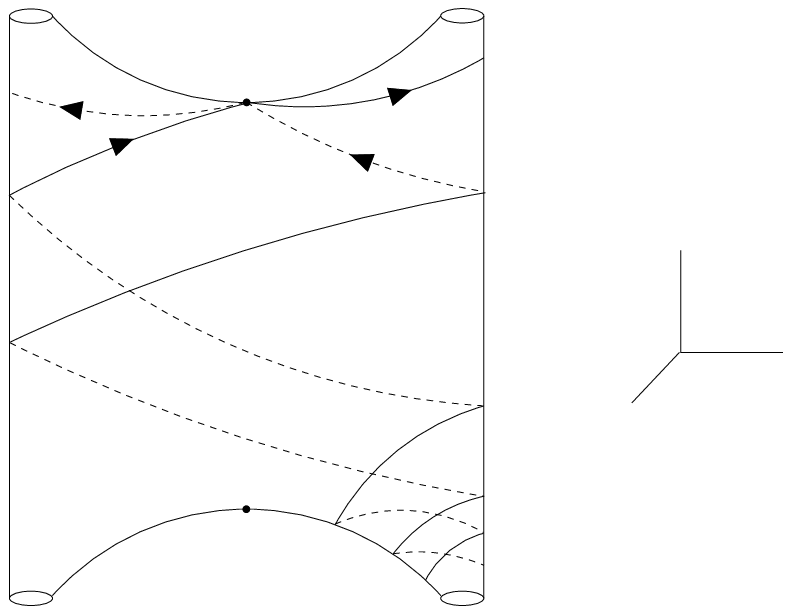}[3.94,2.42]}
\putlabel(1.78,2.06)[b]{$d$}
\putlabel(1.78,.37)[t]{$c$}
\putlabel(3.52,1.42)[b]{$v$}
\putlabel(3.93,1.02)[l]{$s$}
\putlabel(3.3,.82)[tr]{$w$}
\end{picture}
\caption{Collision manifold with stable and unstable branches of $d$ on $\protect\bC$.}\protect\label{fig1}
\end{floatfig}

The function, $W(s)/(1-s^2)$ has exactly one
critical point at $s_c$. Let 
$$
v_c=\sqrt{\frac{W(s_c)}{1-s_c^2}}.
$$
System~\protect\ref{eq:mcg} has two equilibria, both of which lay on the
collision manifold, $\bC$. The $(r, v, s, w)$-coordinates of these two
equilibria are $c=(0,-v_c,s_c,0)$ and $d=(0,v_c,s_,c,0)$ (see
Figure~\protect\ref{fig1})

On the collision manifold, $c$ and $d$ each have a one-dimensional
stable manifold and a one-dimensional unstable manifold. The flow on
the collision manifold is gradient-like with respect to level sets of
$v$ (that is, $dv/d\tau \geq 0$ on $\bC$). A sketch of these stable and
unstable manifolds is given in Figure~\protect\ref{fig1a}.

There is a heteroclinic connection between $c$ and $d$ corresponding
to an ejection-collision orbit (see Figure~\protect\ref{fig1a}). This solution
is homographic, that is, the configuration is constant. So, along the
ejection-collision orbit, the $s$-coordinate is constant, $s_c$. The
forward solution of the orbit limits onto $c$ and the backwards limit
is $d$. Thus the ejection-collision orbit begins and ends at triple
collision and does not pass though a binary collision.

We say any solution which has a forward limit on $c$, ``ends in triple
collision'' and we say any solution which has a backwards limit from
$d$, ``begins in ejection''.  The set of solutions which end in triple
collision forms the two-dimensional stable manifold of $c$, $W^s(c)$
and the set of solutions which begin in ejection forms the
two-dimensional unstable manifold of $d$, $W^u(d)$.  Any orbit in the
intersection of $W^s(c)$ and $W^u(d)$ is called an ejection-collision
solution.

\begin{floatfig}
\begin{picture}(3.94,2.42)
\put(0,-2.42){\loadps{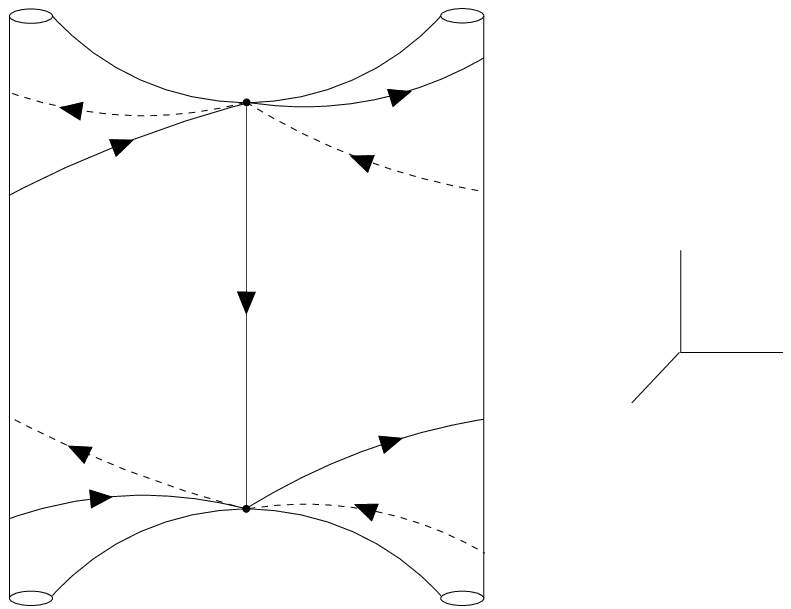}[3.94,2.42]}
\putlabel(1.78,2.06)[b]{$d$}
\putlabel(1.78,.37)[t]{$c$}
\putlabel(3.52,1.42)[b]{$v$}
\putlabel(3.93,1.02)[l]{$s$}
\putlabel(3.3,.82)[tr]{$w$}
\end{picture}
\caption{Ejection-collision orbit with stable and unstable
branches of $c$ and $d$ on $\protect\bC$.} \protect\label{fig1a}
\end{floatfig}

Binary collisions occur at $s=\pm 1$. The energy relation in
Equation~\protect\ref{eq:energy} requires $w=0$ at binary collision. Thus we
must use both the $r$ and $v$ coordinates distinguish one binary
collision from another. Although the System~\protect\ref{eq:mcg} appears to
separate the $r$-coordinate from the others (and indeed does so for
$|s|<1$), we will have to include the $r$-coordinate at binary
collisions.

All solutions, except for the homographic ejection-collision orbit,
pass through at least one binary collision. This makes the pair of
half planes, 
$$
\Gamma=\left\{(r, v, s, w)\:|\: r\geq 0, |s|=1, w=0\right\}
$$ 
the appropriate choice for a
Poincar\'{e} slice. Thanks to work by Mayer and Wang, we know
something about the geometry of the intersection of the stable
manifold for triple collision and this Poincar\'{e} slice, $\Gamma$. 

The motivation for studying how the stable manifold intersect the
Poincar\'{e} slice is given by the following argument. Any point on
the Poincar\'{e} slice has a forward itinerary of binary collisions,
perhaps ending in triple collision. Given a continuous arc on the
Poincar\'{e} slice whose endpoints have different itineraries, there
must be some point on the arc which is also on $W^{s}(c)$. That is,
$W^s(c)$ divides the Poincar\'{e} slice into regions with different
itineraries.

We distinguish the two half-planes of the Poincar\'{e} slice by $L$ for
$s=-1$, and $R$ for $s=1$. The intersections of $W^s(c)$ with $L$ and
$R$ are arcs with two endpoints on $r=0$ or loops with one endpoint on
$r=0$. We label these arcs and loops by their itineraries using $L$'s
and $R$'s for binary collisions and $C$ for triple collision. We use
the subscript $\cd$ to designate the location of the set with a given
itinerary. For example the itinerary $LL\cd RC$ is the set of initial
conditions on $L$ which in forward time pass through $R$ once before
triple collision and in backwards time pass through $L$ again. The
itinerary $LLR\cd C$, however, designates the set of initial conditions
on $R$ which lead directly to triple collision and whose prior two
binary collisions were on $L$. 

Solutions on $\bC$, may lay on the stable manifold for $c$
or for $d$. Such solutions will terminate with the symbol $c$ or $d$
as appropriate. For example, $L\cd Rc$ designates the initial
condition on $\bC$ which traverses to $R$ and then limits onto $c$
without passing though another binary collision. Likewise, $L\cd Ld$
designates the point on $\bC$ which begins on $L$, returns to $L$ and
then limits onto $d$ without passing though another binary
collision. Using this itinerary notation for the intersections of
$W^s(c)$ and the Poincar\'{e} slice, we next summarize Meyer and
Wang's results \protect\cite{MW} in the setting of McGehee coordinates.

There is a unique arc on $L$ with itinerary $L\cd C$ (see
Lemma~\protect\ref{lc-rc}). This arc has two endpoints on $\bC$, one with the
itinerary $L\cd c$ and the other $L\cd d$. Likewise, there is a unique
arc on $R$ with itinerary $R\cd C$. Its endpoints have itineraries
$R\cd c$ and $R\cd d$ (see Figure~\protect\ref{fig2}). We know these arcs are
unique since they must each be homotopic to the ejection-collision
orbit between $c$ and $d$ which passes through no binary collisions.

\begin{center}
\begin{floatfig}
\begin{picture}(3.84,3.31)
\put(0,-3.31){\loadps{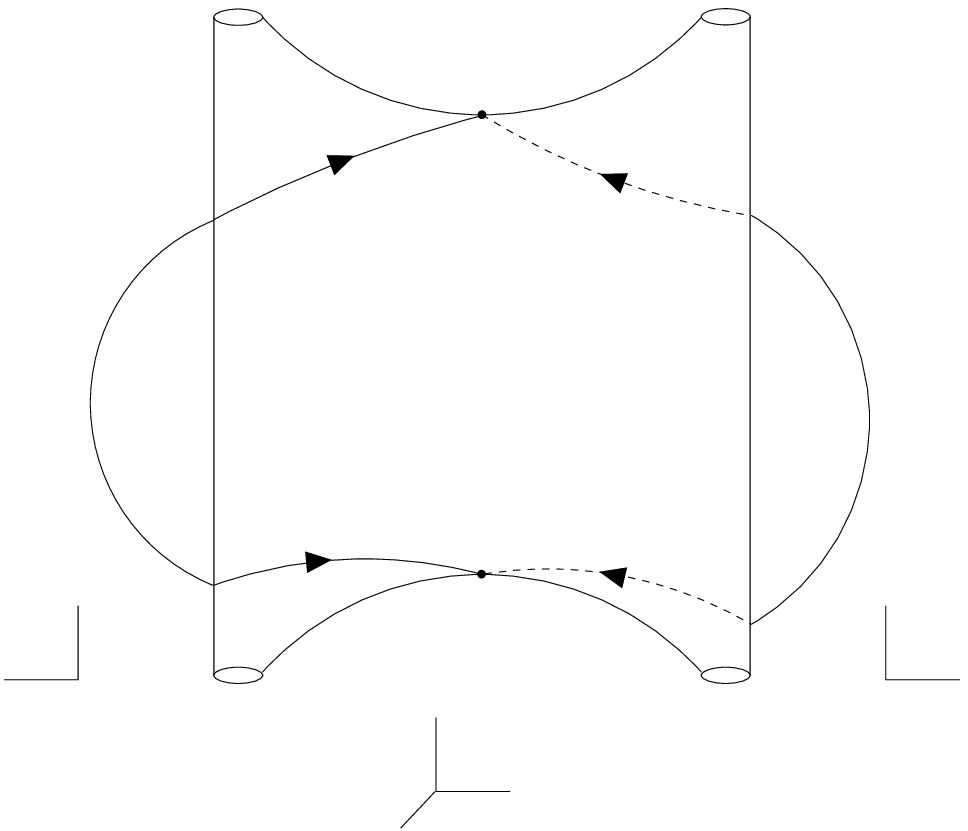}[3.84,3.31]}
\putlabel(0,3)[c]{$L$}
\putlabel(3.84,3)[c]{$R$}
\putlabel(1.92,1)[t]{$c$}
\putlabel(1.92,2.87)[b]{$d$}
\putlabel(1.74,.45)[b]{$v$}
\putlabel(2.05,.15)[l]{$s$}
\putlabel(1.6,0)[tr]{$w$}
\putlabel(3.55,.9)[b]{$v$}
\putlabel(3.84,0.6)[l]{$r$}
\putlabel(.31,0.9)[b]{$v$}
\putlabel(0,.6)[r]{$r$}
\putlabel(.85,.96)[tr]{$L\cd c$}
\putlabel(.85,2.45)[br]{$L\cd d$}
\putlabel(3,.82)[tl]{$R\cd c$}
\putlabel(3,2.47)[bl]{$R\cd d$}
\putlabel(0.37,1.7)[r]{$L\cd C$}
\putlabel(3.5,1.7)[l]{$R\cd C$}
\end{picture}

\caption{First pullback of stable manifold for triple collision.} \protect\label{fig2}
\end{floatfig}
\end{center}

To designate a solution which begins in ejection we use the symbol
$E$. For example, $ELR\cd LC$ designates the initial condition(s) on $R$
which pass though $L$ and limit onto $c$. These initial condition(s)
must also in backwards time pass through $L$ then limit onto $d$
without passing through any other binary collisions.

Intersections of the ejection manifold, $W^u(d)$, with the
Poincar\'{e} slice have a simple relation to intersections of $W^s(c)$
because System~\protect\ref{eq:mcg} is reversible. That is, it has a time
symmetry. Let
\begin{eqnarray*}
\overline{v}&=&-v\\
\overline{w}&=&-w\\
\overline{\tau}&=&-\tau.
\end{eqnarray*}
If $(r, v, s, w)$ is a solution in time $\tau$, then so is $(r,
\overline{v},s, \overline{w})$ in time $\overline{\tau}$.

This means that the arcs $L\cd C$ and $EL\cd$ are reflections of one
another over the $v=0$ axis.  In general the arc or loop with
itinerary $a_0. a_1 a_2...a_k C$ flips about the $v=0$ axis to an arc
or loop with the itinerary $E a_k... a_2 a_1 a_0.$ and vice versa.

With the itinerary notation, we can briefly describe the work of Meyer
and Wang. They generate their main results by focusing on regions
whose itinerary begins $L\cd L$, $L\cd RL$, $L\cd RRL$, and in general
$L\cd R^nL$ (where $R^n$ denotes a string of $n$ $R$'s). These regions
are bounded by pieces of stable manifold whose itinerary is
$LR^nC$. They show that the first return of these regions to $L$
generates a sub-shift on an infinite number of symbols. Since not all
regions are guaranteed to intersect the first returns, the exact
dynamics can not established although it is clear that the dynamics
are rich.

\section{The Stable Manifold of Triple Collision}In this section we characterize pullbacks of the stable manifold for
triple collision, $W^s(c)$, on $L$ and $R$. Recall we are studying the
dynamics of the System~\protect\ref{eq:mcg}. We have chosen a Poincar\'{e} slice,
$\Gamma$, which is transverse to the flow \protect\cite{Mc} made up of two
half-planes, $L$ and $R$, corresponding to left and right binary
collisions. The flow induces a map on $\Gamma$. We denote this
Poincar\'{e} (first return) map by $P$. We now want to show how
the stable manifold for triple collision, $W^s(c)$, intersects
$\Gamma$.

\begin{lemma}[lc-rc]
{The first pull backs of the collision manifold, $L\cd C$ and $R\cd C$, are
smooth arcs with end points $L\cd c$, $L\cd d$ and $R\cd c$, $R\cd d$.}

Since the flow is smooth, $L\cd C$ and $R\cd C$ are deformations of the
ejection-collision orbit between $c$ and $d$. Hence $L\cd C$ and $R\cd C$
are smooth arc with end points on $\bC$. Since the endpoints must
limit to $c$ and $d$, the endpoints of $L\cd C$ are $L\cd C$ and
$L\cd d$. Likewise, the endpoints of $R\cd C$ are $R\cd c$ and $R\cd d$. 
\end{lemma}

\begin{floatfig}
\begin{picture}(2.11,2.63)
\put(0,-2.63){\loadps{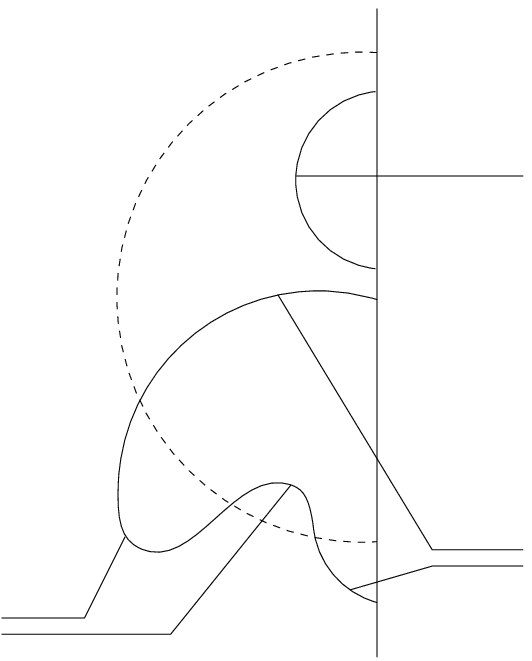}[2.11,2.63]}
\putlabel(2.1,1.93)[l]{Lemma~\ref{case1}}
\putlabel(2.13, .42)[l]{Lemma~\ref{case2}}
\putlabel(-.05,.15)[r]{Lemma~\ref{case3}}
\putlabel(.4,1.6)[r]{$EL\cd$}
\end{picture}

\caption{Typical arcs on $L$ labeled with corresponding Lemmas.}\protect\label{fig3}
\end{floatfig}

We now show how to continue pulling back pieces of stable manifold for
triple collision (see Figure~\protect\ref{fig3}). The following Lemmas are
consistent with Meyer and Wang's results. However, Meyer and Wang
proved them for the case where left binary collisions were not
regularized. Here, we show them for the fully regularized system.

\begin{lemma}[case1]
{If an arc, $\gamma$ in $W^s(c)\cap \Gamma$ has endpoints $x$ and $y$
on $\bC$ and the arc does not intersect $EL\cd$ or $ER\cd$ then the
pullback of $\gamma$, is an arc on $\Gamma$ with endpoints $P^{-1}(x)$
and $P^{-1}(y)$ on $\bC$.}  

Clearly, $x$ pulls back to $P^{-1}(x)$. Points on $\gamma$ near $x$
pull back near $P^{-1}(x)$.  The only obstruction to pulling back a
point on $\gamma$ to $\Gamma$ is if the point pulls back to
ejection. Since $\gamma$ has no such obstruction, all of $\gamma$
pulls back to $\Gamma$ as described.

\end{lemma}

\begin{fullfloatfig}
\begin{picture}(1.06,2.63)
\put(0,-2.63){\loadps{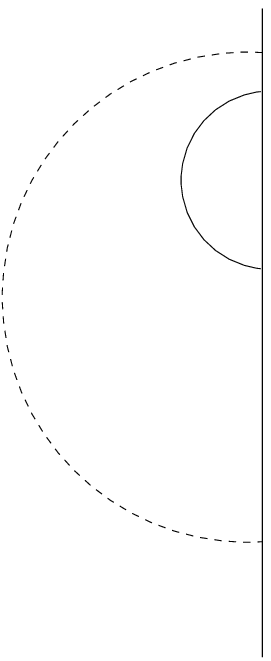}[1.06,2.63]}
\putlabel(0,1.6)[r]{$EL\cd$}
\putlabel(.7,1.9)[r]{$\gamma$}
\putlabel(1.13,2.27)[l]{$x$}
\putlabel(1.13,1.57)[l]{$y$}
\putlabel(0,2.5)[c]{$L$}
\end{picture}
\hspace{1.2in}\begin{picture}(0.36,2.63)
\put(0,-2.63){\loadps{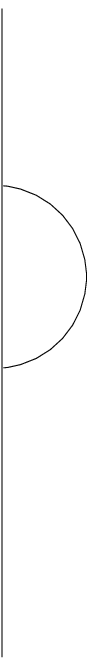}[1.28,2.63]}
\putlabel(-.04,1.18)[r]{$P^{-1}(y)$}
\putlabel(-.04,1.9)[r]{$P^{-1}(x)$}
\putlabel(.42,1.55)[l]{$P^{-1}(\gamma)$}
\putlabel(.3,2.5)[c]{$R$}
\end{picture}

\caption{An arc and its pullback for Lemma~\protect\ref{case1}.}\protect\label{fig4}
\end{fullfloatfig}

\begin{lemma}[case2]
{If an arc in $W^s(c)\cap \Gamma$ has endpoints $x$ and $y$ on $\bC$
and the arc intersects $EL\cd$ or $ER\cd$ then the segment, $\delta$, from
$x$ (or $y$) to the first intersection with $EL\cd$ or $ER\cd$ pulls back
to an arc on $\Gamma$ with endpoints $P^{-1}(x)$ (or $P^{-1}(y)$) and
either $L\cd d$ if the segment pulls back to $L$ or $R\cd d$ if the segment
pulls back to $R$.}

Assume that $\delta$ has an endpoint at $x$. Clearly, $x$ pulls back
to $P^{-1}(x)$. Points on $\delta$ near $x$ pull back near
$P^{-1}(x)$. We continue pulling back points along $\delta$. Points
near the intersection of $\delta$ and $EL\cd $ or $ER\cd $ must pull back
arbitrarily close to the equilibrium, $d$, hence return to $\Gamma$
arbitrarily close to one branch of the stable manifold for $d$ on
$\bC$. Thus points on $\delta$ near the intersection of $\delta$ and
$EL\cd$ or $ER\cd$ must pull back to $\Gamma$ arbitrarily close to either
$L\cd d$ or $R\cd d$. By continuity, points near the intersection of
$\delta$ and $EL\cd $ or $ER\cd $ must pull back to the same half-plane, $L$
or $R$, as the rest of $\delta$.

\end{lemma}

\begin{fullfloatfig}
\begin{picture}(.75,1.83)
\put(0,-1.83){\loadps{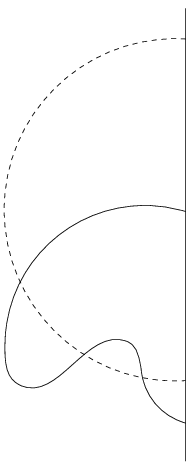}[.75,1.83]}
\putlabel(0,1.05)[r]{$EL\cd$}
\putlabel(.54,1.05)[b]{$\delta_1$}
\putlabel(.7,.2)[tr]{$\delta_2$}
\putlabel(.81,1)[l]{$x$}
\putlabel(.81,.155)[l]{$y$}
\putlabel(0,1.68)[c]{$L$}
\end{picture}
\hspace{1.2in}\begin{picture}(1.75,1.78)
\put(0,-1.78){\loadps{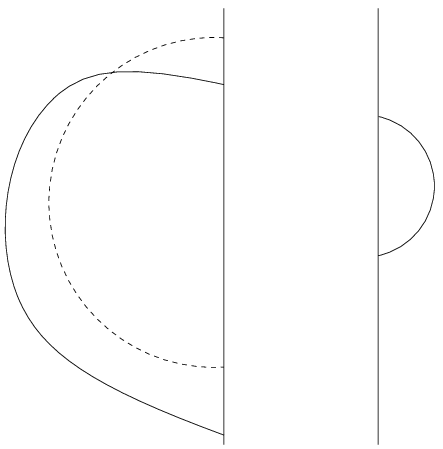}[1.75,1.78]}
\putlabel(.24,1)[l]{$EL\cd$}
\putlabel(1.55,.78)[tl]{$P^{-1}(x)$}
\putlabel(.95,.047)[l]{$P^{-1}(y)$}
\putlabel(1.8,1.05)[l]{$P^{-1}(\delta_1)$}
\putlabel(0,1.05)[r]{$P^{-1}(\delta_2)$}
\putlabel(.95,1.45)[l]{$L\cd d$}
\putlabel(1.55,1.32)[bl]{$R\cd d$}
\putlabel(0,1.68)[c]{$L$}
\putlabel(1.76,1.68)[c]{$R$}
\end{picture}

\caption{Two arc segments and their pullbacks for Lemma~\protect\ref{case2}.}\protect\label{fig5}
\end{fullfloatfig}

\begin{lemma}[case3]
{If a segment of an arc, $\delta$, in $W^s(c)\cap \Gamma$ has both
endpoints on $EL\cd $ or $ER\cd $ and has no other intersections with $EL\cd$
or $ER\cd$, then the segment pulls back to a loop on $\Gamma$ with both
endpoints at either $L\cd d$ if the segment pulls back to $L$ or $R\cd d$ if
the segment pulls back to $R$.}

Since $\delta$ intersects $EL\cd$ or $ER\cd$ only at its endpoints, all of
$\delta$ pulls back to either $L$ or $R$. Points on $\delta$ near the
intersections of $\delta$ and $EL\cd$ or $ER\cd$ must pull back
arbitrarily close to the equilibrium, $d$, hence return to $\Gamma$
arbitrarily close to one branch of the stable manifold for $d$ on
$\bC$. Thus points on $\delta$ near the intersections of $\delta$ and
$EL\cd$ or $ER\cd$ must pull back to $\Gamma$ arbitrarily close to either
$L\cd d$ if the rest of $\delta$ pulls back to $L$ or $R\cd d$ if the rest
of $\delta$ pulls back to $R$.

Pullbacks of loops are covered in earlier Lemmas since we may let
$x=y$ in Lemma~\protect\ref{case1}.
\end{lemma}

\begin{fullfloatfig}
\begin{picture}(0.75,1.83)
\put(0,-1.83){\loadps{fig5a.eps}[0.75,1.83]}
\putlabel(0,1.05)[r]{$EL\cd$}
\putlabel(0,.39)[r]{$\delta_1$}
\putlabel(.47,.51)[b]{$\delta_2$}
\putlabel(0,1.68)[c]{$L$}
\end{picture}
\hspace{1.2in}\begin{picture}(1.26,1.78)
\put(0,-1.78){\loadps{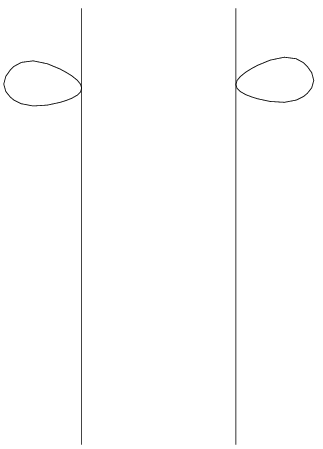}[1.26,1.78]}
\putlabel(0,1.44)[r]{$P^{-1}(\delta_1)$}
\putlabel(1.24,1.44)[l]{$P^{-1}(\delta_2)$}
\putlabel(.34,1.44)[l]{$L\cd d$}
\putlabel(.89,1.44)[r]{$R\cd d$}
\putlabel(0,1.68)[c]{$L$}
\putlabel(1.25,1.68)[c]{$R$}
\end{picture}

\caption{Two arc segments and their pullbacks for Lemma~\protect\ref{case3}.}\protect\label{fig6}
\end{fullfloatfig}

By the lemmas above, we know that the endpoints of arcs and loops in
$W^s(c) \cap \Gamma$ are on either the stable manifold of $c$ or $d$
on $\bC$. By McGehee, we know that the stable manifold for $c$ has two
branches on $\bC$. The stable manifold for $d$ has two branches. In
backwards time, one branch passes first through $L$ and then
alternates between $R$ and $L$, eventually going down one of the
legs. The other branch of $W^s(d)$ on $\bC$ in backwards time passes
first through $R$ and then alternates between $L$ and $R$, eventually
going down one of the legs (see Figure~\protect\ref{fig1}).

One branch of $W^s(d)$ on $\bC$ has an itinerary $...LLLRLRL...Ld$. 
We denote the length of the alternating part of this itinerary
$\lc$. Likewise, we denote the length of the alternating part of the
branch of $W^s(d)$ ending $...Rd$ by $\rc$. For example, a branch of
$W^s(d)$ on $\bC$ with itinerary $...LLLRLRLd$ yields $\lc=5$.

Since the branches of $W^s(d)$ can not intersect we have the relation
\begin{eqnarray}
|\lc -\rc| \leq 1.
\end{eqnarray}
If $\lc=\rc$ then either each branch goes down a different leg of
$\bC$ or there is a heteroclinic connections between $c$ and $d$ on
$\bC$, i.e. $W^s(d)=W^u(c)$. Otherwise $\lc$ and $\rc$ differ by one
and both branches must go down the same leg.

McGehee shows that the $v$-coordinate along flow on $\bC$ is
non-decreasing. Meyer and Wang show that the stable branches for $d$
on $\bC$ must end with $RLd$ or $LRd$. That is, $$\lc, \rc \geq
2.$$ These two facts have the following geometric consequence.

\begin{lemma}[lcandelcross]
{The stable manifold of $c$ with itinerary $L\cd C$ and the unstable
manifold with itinerary $EL\cd$ intersect on $L$. Likewise, $R\cd C$ and
$ER\cd $ intersect on $R$.}

We only look at the case of the intersection of $L\cd C$ and $EL\cd$ on $L$
since the other case is similar.

Let $v(x)$ be the $v$-coordinate of a point on $\bC$. Since $v$ is
non-decreasing along solutions on $\bC$, we know from the definition
of stable and unstable manifolds that the following inequality holds
\begin{eqnarray*}
v(L\cd C)\leq v(c) < 0 < v(d) \leq v(dL\cd)
\end{eqnarray*}

If $v(L\cd d) \leq v(cL\cd)$, then the inverse image of $L\cd d$ under the
Poincar\'{e} map must be on $L$ since the pullback of $L\cd d$ can not
intersect the trajectory between $c$ and $cL\cd$. This means that one
branch of the stable manifold for $d$ on $\bC$ would have the
itinerary $LLd$ which contradicts Meyer and Wang's result. Thus
$v(cL\cd) \leq v(L\cd d)$. Since $v(cL\cd)+v(L\cd d)=0$ we have the inequality
\begin{eqnarray}
v(L\cd C)\leq v(c) \leq v(cL\cd) \leq 0 \leq  v(L\cd d) \leq v(d) \leq v(dL\cd)
\protect\label{inequal}
\end{eqnarray}

The continuous arc $L\cd C$ has endpoints at $L\cd c$ and $L\cd d$. The
continuous arc $EL\cd$ has endpoints at $dL\cd$ and $cL\cd$. By the
inequality in Equation~\protect\ref{inequal}, $L\cd C$ and $EL\cd$ must intersect
at least once.  
\end{lemma}

The segment of $L\cd C$ from $L\cd d$ to the first intersection of $EL\cd$ by
Lemma~\protect\ref{case2} must pullback under $P^{-1}$ to an arc on $R$ with
endpoints at $R\cd Ld$ and $R\cd d$. If this arc does not intersect $ER\cd$
then we continue pulling back this segment until there is a first
intersection with either $EL\cd$ or $ER\cd$. Denote by $\li$ the number of
pullbacks required for the first intersection. Likewise, let $\ri$
denote the number of pullbacks required for the segment of $R\cd C$ from
$R\cd d$ to the first intersection with $ER\cd$ to next intersect $EL\cd$ or
$ER\cd$.

\section{Main Results}The Lemmas of the previous section establish that $W^s(c) \cap \Gamma$
is made up smooth arcs and loops whose endpoints are on the stable
manifolds of $c$ and $d$ on $\bC$. Thus $W^s(c) \cap \Gamma$ generates
an infinite number of regions on $\Gamma$. However, only a finite
number of regions are needed to determine or bound the symbolic
dynamics. 

The values of $\lc, \rc, \li$, and $\ri$ characterize the global
dynamics. The following theorems connect $\lc, \rc, \li$, and $\ri$
with the global dynamics.

Before proceeding, we review the definition of a sofic system. A
finite directed graph whose arrows are labeled and the labels may be
used more than once is called a {\it sofic} system. The dual to a
sofic system, then, is a finite directed graph whose vertices are
labeled and whose labels may be used more than once. Sofic systems are
a generalization of sub-shifts of finite type. See \protect\cite{W} and \protect\cite{LM} for
more information on sofic systems and their role in dynamical
systems. 

\begin{floatfig}
\begin{picture}(3.76,3.65)
\put(0,-3.65){\loadps{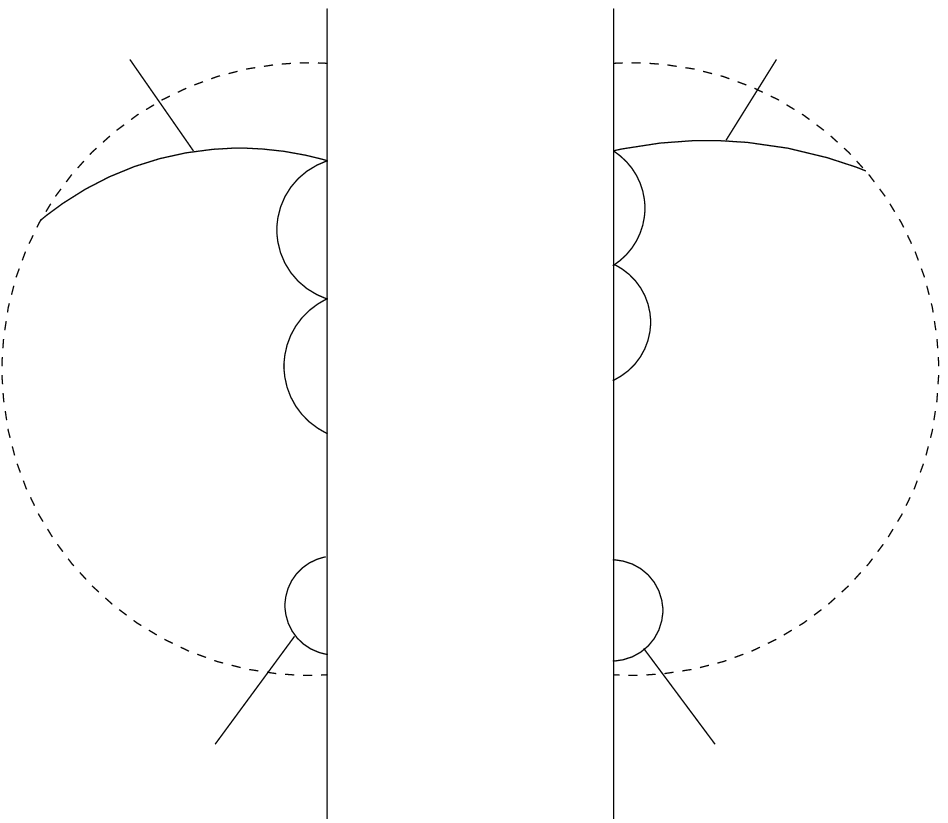}[3.76,3.65]}
\putlabel(.85,.3)[tr]{$P^{-\li+1}(PRS)$}
\putlabel(1.27,1.4)[r]{$\vdots$}
\putlabel(2.5,1.4)[l]{$\vdots$}
\putlabel(.5,3.05)[br]{$PLS$}
\putlabel(3.1,3.1)[b]{$PRS$}
\putlabel(3.76,1.8)[l]{$ER\cd$}
\putlabel(0,1.8)[r]{$EL\cd$}
\putlabel(1.35,2.65)[l]{$L\cd d$}
\putlabel(2.4,2.7)[r]{$R\cd d$}
\putlabel(1.1,2.37)[r]{$P^{-1}(PRS)$}
\putlabel(2.58,2.45)[l]{$P^{-1}(PLS)$}
\putlabel(2.6,2)[l]{$P^{-2}(PRS)$}
\putlabel(1.13,1.8)[r]{$P^{-2}(PLS)$}
\putlabel(2.85,.3)[tl]{$P^{-\li+1}(PLS)$}
\putlabel(1.35,2.05)[l]{$P^{-1}(R\cd d)$}
\putlabel(2.4,2.28)[r]{$P^{-1}(L\cd d)$}
\end{picture}
\caption{Dividing $L$ and $R$ each into $\protect\li+2$ regions. Diagram
labeled for $\protect\li$ even.}  \protect\label{fig7}
\end{floatfig}

\begin{theorem}[main1]
{If for three given masses $\lc=\rc=\li+1=\ri+1$ then the set of
allowed itineraries for the system is given by the dual of a sofic
system (see Figures~\protect\ref{fig9a} and~\protect\ref{fig9b}).}

The arc $EL\cd$ divides $L$ into two regions. We call the bounded region
the {\it inside} of $EL\cd$ and the unbounded region the {\it outside}
of $EL\cd$. Likewise, we define the inside of $ER\cd$ and the outside of
$ER\cd$

We next divide the inside of $EL\cd$ into a finite number of regions
via the following procedure. The segment(s) of $L\cd C$ inside $EL\cd$
is (are) the primary left segment(s), denoted $PLS$. Note that the
segment of $L\cd C$ from from $L\cd d$ to the first intersection with
$EL\cd$ is always contained in the the $PLS$. If $PLS$ contains any
other segments, they must be of the form described by
Lemma~\protect\ref{case3}.

By the definition of $\li$, the first $\li-1$ pullbacks of $PLS$ are
(or at least contain) arcs inside $ER\cd$ and $EL\cd$. We call these
arcs $PLS^{-1}, PLS^{-2}..., PLS^{-\li+1}$ successively. Likewise, we
define $PRS$ and its pullbacks, $PRS^{-1}, PRS^{-2}...,
PRS^{-\li+1}$. The arcs, $EL\cd, ER\cd, PLS$ and its pullbacks and
$PRS$ and its pullbacks divide $L$ and $R$ each into $\li+2$ regions.
Figure~\protect\ref{fig7} describes the case where $PLS$ and $PRS$ each only
have one segment. If they had additional segments, loops would be
attached at the pullbacks of $L\cd d$ and/or $R\cd d$.

Points on $\bC$ of the form $P^{-2i}(L\cd d)$ or $P^{-2i+1}(R\cd d)$
($i\geq0$) which are in the closure of the inside of $EL\cd$ are called
{\it inside points} of $EL\cd$. Likewise, points on $\bC$ of the form
$P^{-2i}(R\cd d)$ or $P^{-2i+1}(L\cd d)$ ($i\geq0$) which are in the closure
of the inside of $ER\cd$ are called inside points of $ER\cd$.
Points on $\bC$ of the form $P^{-i}(R\cd d)$ or $P^{-i}(L\cd d)$ which are
not inside points are called {\it outside points}.

We define a {\it vertical strip} as a simply connected subset of a region
whose intersection with the boundary of the region is made up of one
of the following:
\begin{enumerate}
\ii at least one inside point and a segment of $EL\cd$ or $ER\cd$
\ii at least one outside point and a segment of $EL\cd$ or $ER\cd$
\ii at least two inside points
\ii at least two outside points
\end{enumerate}

The pullback of an inside point is either an inside point or an
outside point. The pullback of an outside point is an outside
point. The pullback of a region arbitrarily near $EL\cd$ or $ER\cd$ is
a region either arbitrarily near $L\cd d$ or $R\cd d$. So the pull
back of a vertical strip at least contains vertical strips. If the
pullback of a vertical strip crosses over $EL\cd$ or $ER\cd$, two
segments of the pullback are vertical strips while the rest of the
segments only have ends on $EL\cd$ or $ER\cd$, and hence are not
vertical (see Figure~\protect\ref{fig8}).

Although we can make a directed graph showing how regions map to one
another under $P^{-1}$, we must also show that the pull back of a
region's vertical strip contain vertical strips in the pullback of the
region. 

\begin{floatfig}
\begin{picture}(3.97,3.27)
\put(0,-3.27){\loadps{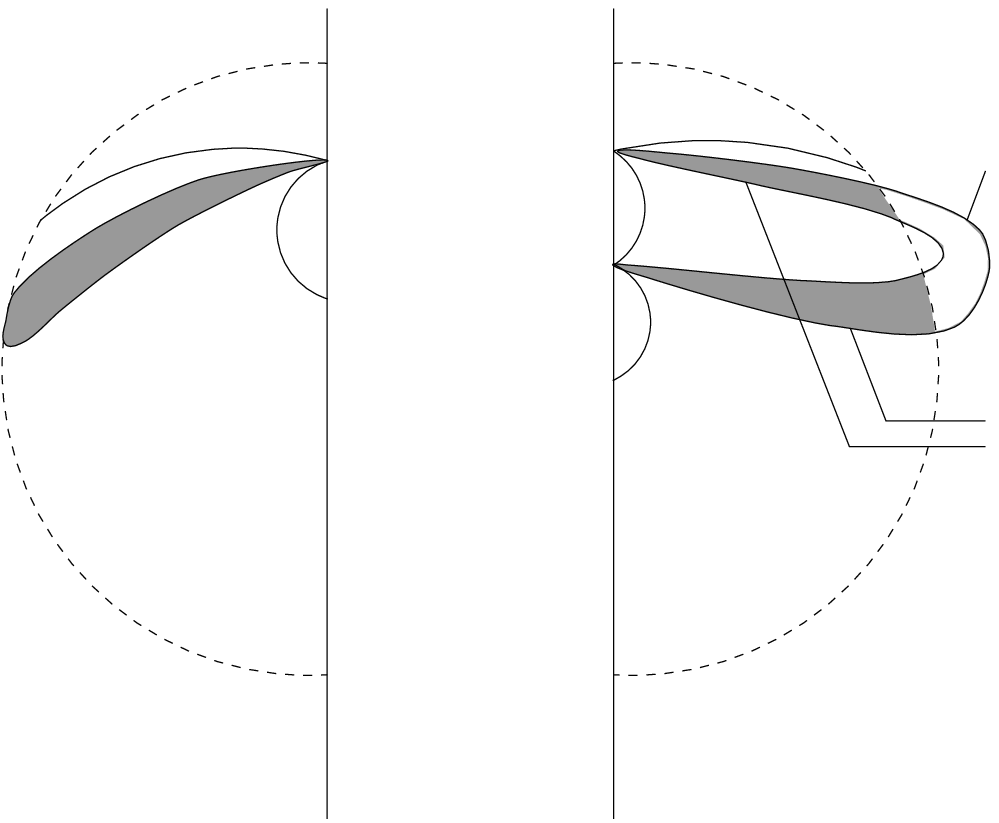}[3.97,3.27]}
\putlabel(2.5,1)[l]{$ER\cd$}
\putlabel(.3,1)[r]{$EL\cd$}
\putlabel(0,3.1)[c]{$L$}
\putlabel(3.76,3.1)[c]{$R$}
\putlabel(1.35,2.65)[l]{$L\cd d$}
\putlabel(2.4,2.7)[r]{$R\cd d$}
\putlabel(3.9,2.64)[bl]{not vertical}
\putlabel(3.9,1.55)[l]{\begin{tabular}{l} 
			pullback has\\ two vertical\\ strips
			\end{tabular}}
\end{picture}

\caption{Shaded vertical strip on $L$ and its pullback on $R$.} \protect\label{fig8}
\end{floatfig}

We next describe the pullbacks of various cases of vertical
strips. Our goal is to show that the division of $L$ and $R$ described
above is a Markov partition of $\Gamma$ for $P^{-1}$ (hence Markov for
$P$ by reversability). 

We begin with a vertical strip inside $EL\cd$ and above the arc of
$PLS$ from $L\cd d$ to the first intersection with $EL\cd$. This
vertical strip pulls back to a vertical strip between $R\cd d$ and
$RL\cd d$ inside $PLS^{-1}$. Further pulls backs generate vertical
strips in regions bounded $PLS^{-i}$ for $1\leq i \leq\li-1$.

A vertical strip inside $PLS^{\li+1}$ contains at least two vertical
strips. One strip includes a segment of $EL\cd$ or $ER\cd$ and the
outside point, $P^{-\li}(L\cd d)$. The other strip includes $EL$ and
the inside point $P^{-\li+1}(R\cd d)$. This second strip, though
inside $EL\cd$ or $ER\cd$ is also outside the regions $PLS^{-i}$ or
$PRS^{-i}$ for $1\leq i \leq\li-1$.

A vertical strip from $EL\cd$ or $ER\cd$ to any outside point has a
pullback which includes a vertical strip from inside $EL\cd$ or
$ER\cd$ and above either the arc of $PLS$ from $L\cd d$ to the first
intersection with $EL\cd$ or the arc of $PRS$ from $R\cd d$ to the
first intersection with $ER\cd$.

Like a vertical strip inside $PLS^{\li+1}$, a vertical strip from a
segment of $EL\cd$ or $ER\cd$ and with $P^{-\li+1}(L\cd d$ or
$P^{-\li+1}(R\cd d)$ contains at least two vertical
strips. One strip includes a segment of $EL\cd$ or $ER\cd$ and the
outside point, $P^{-\li}(L\cd d)$. The other strip includes $EL$ and
the inside point $P^{-\li+1}(R\cd d)$. This second strip, though
inside $EL\cd$ or $ER\cd$ is also outside the regions $PLS^{-i}$ or
$PRS^{-i}$ for $1\leq i \leq\li-1$.

The above argument can be repeated beginning with a vertical strip
inside $ER\cd$ and above the arc of $PRS$ from $R\cd d$ to the first
intersection with $ER\cd$.

\begin{floatfig}
\vspace{.1in}
\begin{picture}(3.00,1.31)
\put(0,-1.31){\loadps{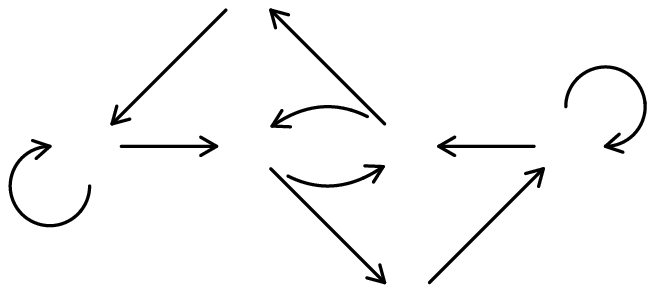}[3.00,1.31]}
\putlabel(.70,0.66)[c]{$R$}           
\putlabel(1.34,0.66)[c]{$L$}           
\putlabel(1.97,0.66)[c]{$R$}           
\putlabel(2.61,0.66)[c]{$L$}           
\putlabel(1.39,1.35)[c]{$A_L^{\li-1}$} 
\putlabel(2.02,-.02)[c]{$A_R^{\li-1}$}  
\end{picture}
\caption{Directed graph for
Theorem~\protect\ref{main1} where $\protect\li$ is even.}  \protect\label{fig9a}
\end{floatfig}

Since vertical strips in a region pull back to vertical strips in the
pullback of the region, the partition of $L$ and $R$ given above is
Markov. If we label the regions, then there is a subshift of finite
type which semi-conjugate to $P^{-1}$ on $\Gamma$. By reversability,
we can switch the direction of the arrows and get a subshift of finite
type which is semi-conjugate to $P$ on $\Gamma$.

Replacing region names with ``$L$'' for regions on $L$ and ``$R$'' for
regions on $R$ yields a directed graph on two symbols (the dual of a
sofic system). This directed graph in general can be simplified. The
simplified graph describes all itineraries of binary collisions.

To see this last claim, any initial condition on $L$ or $R$ is in some
region as defined above. Thus the point maps about the regions
according to some path in the sub-shift, so its itinerary is described
by the graph.

Conversely, given any itinerary described by the graph, there is at
least one path in the sub-shift which accomplishes the given
itinerary. Choose a closed vertical strip in the first region for such
a path. Since the pullback of vertical strips includes contains
vertical strips we can generate a nested sequence of closed vertical
strips which obtain arbitrarily many terms in the desired
sequence. The infinite intersection of these closed and nested
vertical strips is non-empty, thus guaranteeing at least one point
which achieves the given itinerary.
\end{theorem}

Carrying out the procedure outlined in Theorem~\protect\ref{main1} we see that
for the case $\lc=\rc=\li+1=\ri+1$, the generated dual to a sofic
system is given by Figures~\protect\ref{fig9a}~and~\protect\ref{fig9b}. Note that the
symbol $A_L^{\li-1}$ means an alternating sequence of $L$'s and $R$'s
beginning with $L$ of length $\li-1$. Likewise, the
symbol $A_R^{\li-1}$ means an alternating sequence of $L$'s and $R$'s
beginning with $R$ of length $\li-1$.

\begin{floatfig}
\vspace{.1in}
\begin{picture}(3.00,1.31)
\put(0,-1.31){\loadps{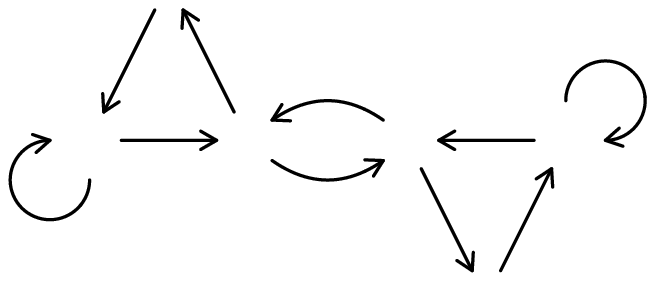}[3.00,1.31]}
\putlabel(.70,0.66)[c]{$R$}
\putlabel(1.34,0.66)[c]{$L$}
\putlabel(1.97,0.66)[c]{$R$}
\putlabel(2.61,0.66)[c]{$L$}
\putlabel(1.07,1.35)[c]{$A_R^{\li-1}$}
\putlabel(2.34,-0.02)[c]{$A_L^{\li-1}$}
\end{picture}
\caption{Directed graph for
Theorem~\protect\ref{main1} where $\protect\li$ is odd.} \protect\label{fig9b}
\end{floatfig}

Note that for $\lc=\rc=\li+1=\ri+1=2$ the dual to a sofic system
reduces to the full shift on two symbols.

\begin{theorem}[main2]
{If for three given masses $\lc=\rc+1=\ri+2\leq \li$ or
$\rc=\lc+1=\li+2\leq \ri$ then the set of itineraries for the system
is bounded between two duals of sofic systems (see
Figures~\protect\ref{fig10a} to~\protect\ref{fig10d}).}

We will assume the first case, $\lc=\rc+1=\ri+2\leq \li$, since the
argument for the second case is similar. 

Using the notation from Theorem~\protect\ref{main1} we divide the inside of
$EL\cd$ and $ER\cd$ into a finite number of regions bounded by
$PLS^{-i}$ for $0\leq i\leq \ri$ and $PRS^{-i}$ for $0\leq i\leq
\ri-1$. We also make regions outside $ER\cd$ or $EL\cd$
bounded by $PLS^{-i}$ for $\ri+1 \leq i \leq \li-1$.

By the argument in Theorem~\protect\ref{main1}, a region's vertical strip pulls
back to vertical strips in the region's pullbacks with one
exception. A vertical strip inside $PLS^{-\li+1}$ will not pull back
to a vertical strip inside $EL\cd$ or $ER\cd$ even though the region
bounded by $PLS^{-\li+1}$ and $\bC$ must cross either $EL\cd$ or
$ER\cd$. This means that a directed graph showing how regions map under
$P^{-1}$ will not be Markov. However a point on $\Gamma$ must move
from region to region according to the the directed graph so the
directed graph contains all allowed itineraries of regions even
though some itineraries may not be achieved. 

So we begin with the directed graph showing how regions map under
$P^{-1}$. By reversability, switching the direction of the arrows we
have a directed graph showing how regions map under $P$. Replacing
region names with ``$L$'' for regions on $L$ and ``$R$'' for regions on $R$
yields a directed graph on two symbols (the dual of a sofic
system). This directed graph in general can be simplified. The
simplified graph contains all itineraries of binary collisions.

To bound the set of allowed itineraries from below, we assume that
pullbacks of $PLS$ never intersect $EL\cd$ or $ER\cd$. That is we assume
$\ri=+\infty$. Although this is not the case, the condition guarantees
a region's vertical strips pull back to vertical strips in the
pullback of the region. That is, the directed graph showing how
regions map under $P^{-1}$ with the condition that $\ri=+\infty$ is
Markov, hence the sub-shift contains orbits which are guaranteed to
occur in the actual system. 

\begin{fullfloatfig}
\begin{picture}(3.00,1.87)
\put(0,-1.87){\loadps{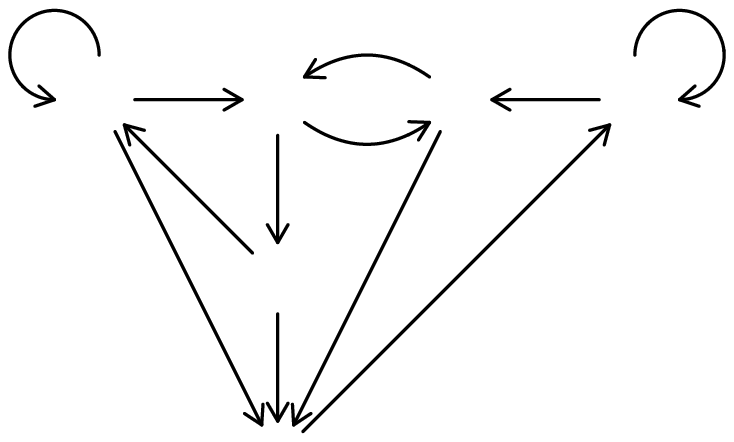}[3.00,1.87]}
\putlabel(0.42,1.47)[c]{$R$}
\putlabel(1.14,1.47)[c]{$L$}
\putlabel(1.85,1.47)[c]{$R$}
\putlabel(2.57,1.47)[c]{$L$}
\putlabel(1.19,0.75)[c]{$A_R^{\lc-1}$}
\putlabel(1.19,.04)[c]{$R^{\li-1}$}
\end{picture}
\caption{Directed graph for Theorem~\protect\ref{main2} containing allowed
sequences for $\protect\lc$ even.}\protect\label{fig10a}
\begin{picture}(3.00,1.87)
\put(0,-1.87){\loadps{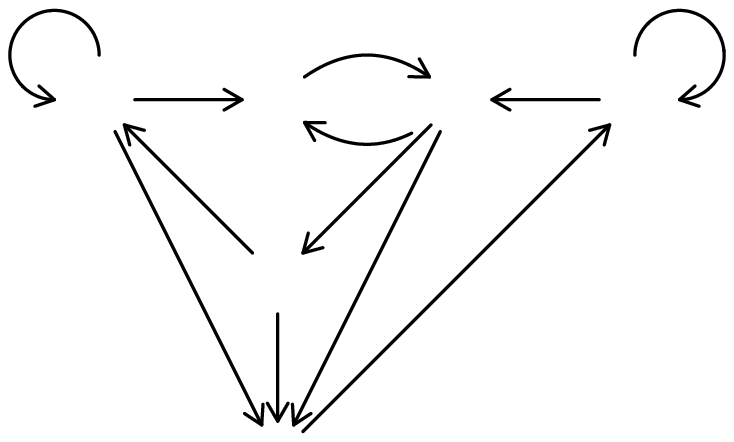}[3.00,1.87]}
\putlabel(0.42,1.47)[c]{$R$}
\putlabel(1.14,1.47)[c]{$L$}
\putlabel(1.85,1.47)[c]{$R$}
\putlabel(2.57,1.47)[c]{$L$}
\putlabel(1.19,0.75)[c]{$A_L^{\lc-1}$}
\putlabel(1.19,.04)[c]{$R^{\li-1}$}
\end{picture}
\caption{Directed graph for Theorem~\protect\ref{main2} containing allowed
sequences for $\protect\lc$ odd.}\protect\label{fig10b}
\end{fullfloatfig}

We begin with the directed graph showing regions map under $P^{-1}$
with the condition that $\ri=+\infty$. By reversability, switching the
direction of the arrows we have a directed graph showing how regions
map under $P$. Replacing region names with ``$L$'' for regions on $L$ and
``$R$'' for regions on $R$ yields a directed graph on two symbols (the
dual of a sofic system). This directed graph in general can be
simplified. The simplified graph contains all guaranteed itineraries
of binary collisions.
\end{theorem}

Carrying out the procedure outlined in Theorem~\protect\ref{main2} we see that
for the case $\lc=\rc+1=\ri+2\leq \li$, the generated duals to sofic
systems is given by Figures~\protect\ref{fig10a}~to~\protect\ref{fig10d}. For the case
$\rc=\lc+1=\li+2\leq \ri$, the generated duals to sofic systems is
given by Figures~\protect\ref{fig10a}~to~\protect\ref{fig10d} after exchanging $L$'s
and $R$'s.

\begin{theorem}[main3]
{If for three given masses, if $\lc, \rc, \li$ and $\ri$ do not meet
the criteria of Theorem~\protect\ref{main1} or Theorem~\protect\ref{main2} then the
set of itineraries for the system is bounded between two duals of
sofic systems.}

If $\lc, \rc, \li$ and $\ri$ do not meet the criteria of
Theorem~\protect\ref{main1} or Theorem~\protect\ref{main2} then (using the notation of
Theorem~\protect\ref{main1}) one of the pull backs of $PLS$ or $PRS$ must
cross $EL\cd$ or $ER\cd$ when both of endpoints of that arc are either
inside $EL\cd$ or $ER\cd$. This means that a vertical strip which
landed in this region will pass to a non-vertical strip outside
$EL\cd$ or $ER\cd$, even though the region pulls back outside $EL\cd$
or $ER\cd$. This means a Markov partition is not possible for a
directed graph on any finite set of regions similar to
Theorem~\protect\ref{main1} or Theorem~\protect\ref{main2}.

If $\lc=\rc$ then the directed graph from Theorem~\protect\ref{main1} 
for
these values of $\lc$ and $\rc$ must contain guaranteed dynamics for
our system since the premature intersection described above only adds
to the possible set of itineraries. 

Likewise, if $\lc \neq \rc$ then the associated directed graph from
Theorem~\protect\ref{main2} which contains guaranteed sequences also contains
guaranteed sequences for our system since the premature intersection
described above only adds to the possible set of itineraries.

\begin{fullfloatfig}
\begin{picture}(3.00,1.50)
\put(0,-1.5){\loadps{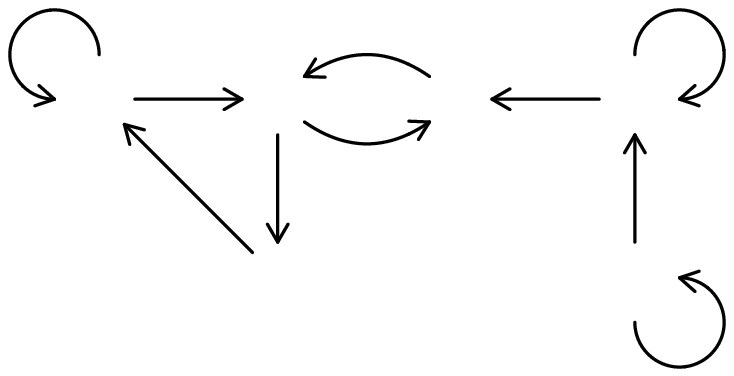}[3.00,1.50]}
\putlabel(0.42,1.10)[c]{$R$}
\putlabel(1.14,1.10)[c]{$L$}
\putlabel(1.85,1.10)[c]{$R$}
\putlabel(2.57,1.10)[c]{$L$}
\putlabel(1.04,0.39)[c]{$A_R^{\lc-1}$}
\putlabel(2.57,0.39)[c]{$R$}
\end{picture}
\caption{Directed graph for Theorem~\protect\ref{main2} containing guaranteed
sequences for $\protect\lc$ even.}\protect\label{fig10c}
\begin{picture}(3.00,1.50)
\put(0,-1.5){\loadps{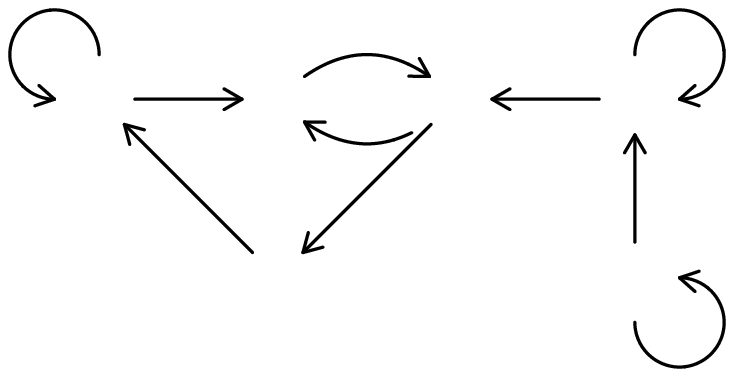}[3.00,1.50]}
\putlabel(0.42,1.10)[c]{$R$}
\putlabel(1.14,1.10)[c]{$L$}
\putlabel(1.85,1.10)[c]{$R$}
\putlabel(2.57,1.10)[c]{$L$}
\putlabel(1.24,0.39)[c]{$A_L^{\lc-1}$}
\putlabel(2.57,0.39)[c]{$R$}
\end{picture}
\caption{Directed graph for Theorem~\protect\ref{main2} containing guaranteed
sequences for $\protect\lc$ odd.}\protect\label{fig10d}
\end{fullfloatfig}

To get an upper bound for our system we divide $L$ and $R$ into a
finite number of regions bounded by $EL\cd$, $ER\cd$, $PLS^{-i}$ for
$0\leq i \leq \li-1$ and $PRS^{-i}$ for $0\leq i \leq \ri-1$. Each of these
pullbacks are single arcs by the conditions on $\li$ and $\ri$. The
directed graph describing how these regions map under $P^{-1}$ must
contain all allowed itineraries of regions even though some
itineraries may not be achieved.

So we begin with the directed graph showing how regions map under
$P^{-1}$. By reversability, switching the direction of the arrows we
have a directed graph showing how regions map under $P$. Replacing
region names with ``$L$'' for regions on $L$ and ``$R$'' for regions on $R$
yields a directed graph on two symbols (the dual of a sofic
system). This directed graph in general can be simplified. The
simplified graph contains all itineraries of binary collisions.
\end{theorem}

Carrying out the procedure outlined in Theorem~\protect\ref{main3} as an
example, we see that for the case $\lc=5, \rc=6, \li=2$ and $\ri=7$,
the generated dual to a sofic system containing all allowed sequences is given
by Figure~\protect\ref{fig11}.

\begin{floatfig}
\begin{picture}(3.00,1.91)
\put(0,-1.91){\loadps{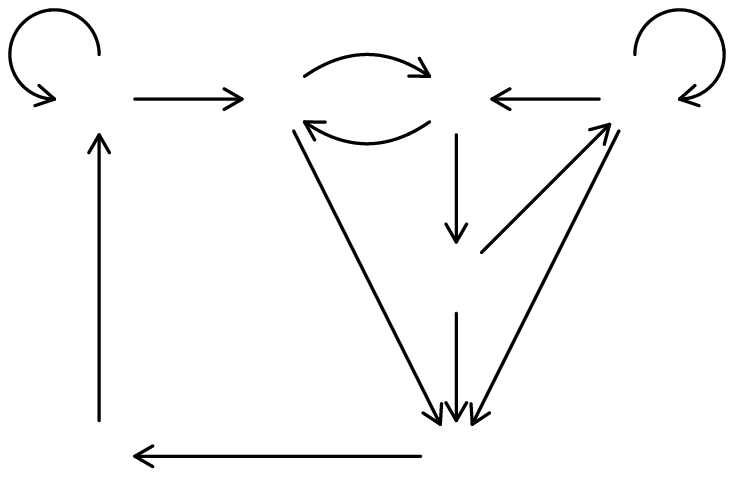}[3.00,1.91]}
\putlabel(0.42,1.51)[c]{$R$}
\putlabel(1.14,1.51)[c]{$L$}
\putlabel(1.85,1.51)[c]{$R$}
\putlabel(2.57,1.51)[c]{$L$}
\putlabel(1.85,0.80)[c]{$L$}
\putlabel(0.34,.08)[c]{$A_L^6$}
\putlabel(1.85,.08)[c]{$L$}
\end{picture}

\caption{Directed graph for Theorem~\protect\ref{main3} containing all allowed sequences for the
case $\protect\lc=5, \protect\rc=6, \protect\li=2$ and $\protect\ri=7$.}\protect\label{fig11}
\end{floatfig}

Note that for $\li=\ri=1$ the upper bound for the dynamics is the full
shift on two symbols, hence a trivial upper bound.

\section{Oscillatory Motion}One application of the Theorems in Section~5 is to the presence of
oscillatory motion in the $N$-body problem. Saari and Xia explored
possible behaviors in the $N$-body problem as $t \rightarrow \infty$
\protect\cite{SX}. One such behavior is oscillatory motion, that is, a mutual
distance coordinate $r$, so that as $t\rightarrow\infty$, the $\limsup
r=\infty$ while the $\liminf r <\infty$. 

For the collinear three-body problem, oscillatory motion requires that
either $m_1$ or $m_3$ takes longer and longer excursions for the other
binary pair, each time returning to the binary pair before its next
excursion. Saari and Xia showed that it is enough to prove the
presence of itineraries of the form
$$
...L^{a_1}...L^{a_2}...L^{a_3}...L^{a_4}... \qquad\mbox{or}\qquad 
...R^{a_1}...R^{a_2}...R^{a_3}...R^{a_4}...
$$
so that $a_i \rightarrow \infty$ as $i\rightarrow\infty$. The
existence of such itineraries guarantees the existence of oscillatory
motion.

Saari and Xia study the return map on the zero momentum set for
$m_3$. To guarantee a return map, they restrict the masses to the case
where triple collision can lead to arbitrarily high velocities. This
is exactly when the stable manifold of $d$ and and the unstable
manifold of $c$ intersect transversely. In this case, the set if
allowed itineraries is given by the full-shift on two symbols (same as
$\lc=\rc=\li+1=\ri+1=2$).

Once they establish the existence of oscillatory motion in the
collinear three-body problem for some sets of masses, they note that
the motion can be extended to the $N$-body problem. They conclude that
for any $N\geq 3$ there exist masses and initial conditions so that
oscillatory motion exists.

From the Theorems in Section~5 it is clear that for all sets of masses
which admit transverse intersections of the stable manifold of $d$ and
and the unstable manifold of $c$, itineraries of the form
$$
...L^{a_1}...L^{a_2}...L^{a_3}...L^{a_4}... \qquad\mbox{or}\qquad 
...R^{a_1}...R^{a_2}...R^{a_3}...R^{a_4}...
$$
so that $a_i \rightarrow \infty$ as $i\rightarrow\infty$ exist. We are
thus led to the same conclusion as a Corollary to our Theorems. 

\statecorollary{For the $N$-body problem, $N\geq 3$, there exist
positive masses and initial conditions so that oscillatory motion
occurs.}

\end{document}